\documentclass[12pt]{article}
\usepackage[dvipdfm]{graphicx}          
\usepackage{tikz}

 \usepackage{latexsym}
 \usepackage{amssymb}
 \usepackage{color}
 \usepackage{graphicx}
\usepackage{epstopdf}
\usepackage{amsmath}    
\usepackage{bm}                     
 \newtheorem{Theorem}{Theorem}
\newtheorem{Definition}{Definition}

\newtheorem{Conjecture}{Conjecture}
\newtheorem{Lemma}{Lemma}
\newtheorem{Corollary}{Corollary}

\newcommand{\A}{{\cal A}}

\newcommand{\vv}{{\bf v}}
\newcommand{\x}{{\bf x}}

\newcommand{\e}{{\bf e}}

\newcommand{\qed}{\nobreak \ifvmode \relax \else
      \ifdim\lastskip<1.5em \hskip-\lastskip
      \hskip1.5em plus0em minus0.5em \fi \nobreak
      \vrule height0.75em width0.5em depth0.25em\fi}

\def \ep{\hbox{ }\hfill$\Box$}

\begin{document}
\title{Positive Semi-Definiteness and Sum-of-Squares Property of Fourth Order Four Dimensional Hankel Tensors}

\author{Yannan Chen\footnote{School of Mathematics and Statistics,
    Zhengzhou University, Zhengzhou, China. E-mail: ynchen@zzu.edu.cn (Y. Chen)
    This author's work was supported by the
    National Natural Science Foundation of China (Grant No. 11401539) and
    the Development Foundation for Excellent Youth Scholars of
    Zhengzhou University (Grant No. 1421315070).}
    \quad Liqun Qi\footnote{Department of Applied Mathematics,
    The Hong Kong Polytechnic University, Hung Hom, Kowloon, Hong Kong.
    E-mail: maqilq@polyu.edu.hk (L. Qi).
    This author's work was partially supported by
    the Hong Kong Research Grant Council (Grant No. PolyU 502111, 501212, 501913 and 15302114).}
    \quad Qun Wang \footnote{Department of Applied Mathematics,
    The Hong Kong Polytechnic University, Hung Hom, Kowloon, Hong Kong.
    Email: wangqun876@gmail.com (Q. Wang).} }

\date{\today} \maketitle

\begin{abstract}
\noindent  
  A positive semi-definite (PSD) tensor which is not a sum-of-squares (SOS) tensor
  is called a PSD non-SOS (PNS) tensor.
  Is there a fourth order four dimensional PNS Hankel tensor?
  Until now, this question is still an open problem.
  Its answer has both theoretical and practical meanings.
  We assume that the generating vector $\vv$ of the Hankel tensor $\A$ is symmetric.
  Under this assumption, we may fix the fifth element $v_4$ of $\vv$ at $1$.
  We show that there are two surfaces $M_0$ and $N_0$ with
  the elements $v_2, v_6, v_1, v_3, v_5$ of $\vv$ as variables,
  such that $M_0 \ge N_0$, $\A$ is SOS if and only if $v_0 \ge M_0$,
  and $\A$ is PSD if and only if $v_0 \ge N_0$, where $v_0$ is the first element of $\vv$.
  If $M_0 = N_0$ for a point $P = (v_2, v_6, v_1, v_3, v_5)^\top$, then
  there are no fourth order four dimensional PNS Hankel tensors with symmetric generating vectors
  for such $v_2, v_6, v_1, v_3, v_5$. Then, we call such a point $P$ PNS-free.
  We show that a $45$-degree planar closed convex cone, a segment,
  a ray and an additional point are PNS-free.
  Numerical tests check various grid points, and find that they are also PNS-free.

\noindent {\bf Key words:}\hspace{2mm} Hankel tensor, generating vector,
  sum of squares, positive semi-definiteness, PNS-free.

\noindent {\bf AMS subject classifications (2010):}\hspace{2mm} 15A18; 15A69
  \vspace{3mm}

\end{abstract}

\section{Introduction}
\hspace{4mm}

In 1888, young Hilbert \cite{Hi} proved that for homogeneous polynomials, only in the following three cases, a positive semi-definite (PSD) polynomial definitely is a sum-of-squares (SOS) polynomial: 1) $m = 2$; 2) $n = 2$; 3) $m=4$ and $n=3$, where $m$ is the degree of the polynomial and $n$ is the number of variables.    Hilbert proved that in all the other possible combinations of $n$ and even $m$, there are PSD non-SOS (PNS) homogeneous polynomials.  The most well-known PNS homogeneous polynomial is the Motzkin polynomial \cite{Mo} with $m=6$ and $n=3$.   Other examples of
PNS homogeneous polynomials was found in \cite{AP, Ch, CL, Re}.

A homogeneous polynomial is uniquely corresponding to a symmetric tensor \cite{Qi}. For a symmetric tensor, $m$ is its order and $n$ is its dimension.   One class of symmetric tensors are Hankel tensors. Hankel tensors arise from signal processing and some other applications \cite{BB, DQW, PDV, Qi15, Xu}.  In \cite{Qi15}, two classes of PSD Hankel tensors were identified.   They are even order strong Hankel tensors and even order complete Hankel tensors.    It was proved in \cite{LQX} that complete Hankel tensors are strong Hankel tensors, and even order strong Hankel tensors are SOS tensors.
It was also shown there that there are SOS Hankel tensors and PSD Hankel tensors, which are not strong Hankel tensors.
Thus, a question was raised in \cite{LQX}:  Are all PSD Hankel tensors SOS tensors?   If there are no PSD non-SOS Hankel tensors, the problem for determining a given even order Hankel tensor is PSD or not can be solved by solving a semi-definite linear programming  problem \cite{LQX, Las, Lau}.

We may call the problem raised by the above question as the Hilbert-Hankel problem,
as in a certain sense, it is the Hilbert problem with the Hankel constraint.

According to Hilbert \cite{Hi, Re}, one case with low values of $m$ and $n$, in which there are PNS homogeneous polynomials, is that $m=6$ and $n=3$.   In \cite{LQW}, the Hilbert-Hankel problem with order six and dimension three was studied.   Four special cases were analyzed.   Thousands of random examples were checked.  No PNS Hankel tensors of order six and dimension three were found in \cite{LQW}.  Theoretically, it is still an open problem whether there are PNS Hankel tensors of order six and dimension three or not.

According to Hilbert \cite{Hi, Re}, another case with low values of $m$ and $n$, in which there are PNS homogeneous polynomials, is that $m=n=4$.   In this paper, we consider Hankel tensors of order four and dimension four.

Let $\vv = (v_0, v_1,\ldots, v_{12})^\top \in \Re^{13}$.   A fourth order four dimensional {\bf Hankel tensor} $\A = (a_{i_1i_2i_3i_4})$ is defined by $$a_{i_1i_2i_3i_4} = v_{i_1+i_2+i_3+i_4-4},$$
for $i_1, i_2, i_3, i_4 = 1, 2, 3, 4$.  The corresponding vector $\vv$ that defines the Hankel tensor $\A$ is called the {\bf generating vector} of $\A$.
For $\x = (x_1, x_2, x_3, x_4)^\top \in \Re^4$, a Hankel tensor $\A$ uniquely defines a Hankel polynomial
\begin{equation} \label{e1}
f(\x) \equiv \A \x^{\otimes 4} = \sum_{i_1, i_2, i_3, i_4=1}^4 a_{i_1i_2i_3i_4}x_{i_1}x_{i_2}x_{i_3}x_{i_4} = \sum_{i_1, i_2, i_3, i_4=1}^4 v_{i_1+i_2+i_3 +i_4-4}x_{i_1}x_{i_2}x_{i_3}x_{i_4}.
\end{equation}
If $f(\x) \ge 0$ for all $\x \in \Re^4$, the Hankel tensor $\A$ is called {\bf positive semi-definite} (PSD).
If $f(\x)$ can be represented as a sum of squares of quadratic homogeneous polynomials,
the Hankel tensor $\A$ is called {\bf sum-of-squares} (SOS).
Clearly, $\A$ is PSD if it is SOS.

In the next section, we present some necessary conditions for the positive semi-definiteness of fourth order four dimensional Hankel tensors.

We may see that the role of $v_j$ is symmetric in $f(\x)$.  In Section 3, we assume that
\begin{equation} \label{e2}
v_j = v_{12-j}
\end{equation}
for $j = 0, \ldots, 5$.  Under this assumption, by the results of Section 2, if $\A$ is PSD, we have $v_0 = v_{12} \ge 0$ and $v_4 = v_8 \ge 0$.    Moreover, if $v_4 = v_8 = 0$ and $\A$ is PSD, $\A$ is SOS.   Thus, we may only consider the case that $v_4 = v_8 > 0$.   Since $\A$ is PSD or SOS or PNS if and only if $\alpha \A$ is PSD or SOS or PNS respectively, where $\alpha$ is an arbitrary positive number, we may simply assume that
\begin{equation} \label{e3}
v_4 = v_8 = 1.
\end{equation}
Next, we show that there is a function $\eta(v_5, v_6)$ such that $\eta(v_5, v_6) \le 1$ if $\A$ is PSD.    We propose that there are two functions $M_0(v_2, v_6, v_1, v_3, v_5) \ge N_0(v_2, v_6, v_1, v_3, v_5)$, defined for $\eta(v_5, v_6) < 1$, such that $\A$ is SOS if and only if $v_0 \ge M_0$, and $\A$ is PSD if and only if $v_0 \ge N_0$. If $M_0 = N_0$ for some $v_2, v_6, v_1, v_3, v_5$, then there are no fourth order four dimensional PNS Hankel tensors for such $v_2, v_6, v_1, v_3, v_5$ under the symmetric assumption (\ref{e2}).   We call such a point $P = (v_2, v_6, v_1, v_3, v_5)^\top \in \Re^5$ a {\bf PNS-free point} of fourth order four dimensional Hankel tensors, or simply a PNS-free point.   We call the set of points in $\Re^5$, satisfying $\eta(v_5, v_6) < 1$, the {\bf effective domain} of fourth order four dimensional Hankel tensors, or simply the effective domain, and denote it by $S$.  We show that if all the points in $S$ are PNS-free, then there are no fourth order four dimensional PNS Hankel tensors with symmetric generating vectors.

\begin{figure}
\begin{center}
\begin{tikzpicture}
    \clip (-5,-2) rectangle (5.7,5.6);
    \draw[step=1cm, brown, very thin] (-4.9,-1) grid (4.9,4.9);
    \draw[->] (-5,0) -- (5.1,0) node[right]{$v_2$};
    \draw[->] (0,-2) -- (0,5.1) node[above]{$v_6$};

    \draw[line width=3pt, red] (-5,0) -- (0.05,0);
    \node[above] at (-2.5,0) {\footnotesize $R$};

    \draw[fill=red, red] (5,1) -- (1,1) -- (5,5);
    \node[right, fill=white] at (2.8,2.1) {\footnotesize PSD=SOS};
    \node[above] at (2.2,1.3) {\footnotesize $C$};

    \draw[fill=gray, gray] (-5,-1/3) -- (-5,-2) -- (5,-2) -- (5,-1/3) -- cycle;
    \node[below, fill=white] at (-2,-.8) {\footnotesize Non PSD};

    \shade[ball color=red] (1,0) circle (3pt);
    \node[below] at (1,0) {\footnotesize $A$};

    \foreach \i in {.8,0.79,...,0} {\shade[ball color=red!50] (1+\i*.574,1-\i*.819) circle (1.5pt);}
    \foreach \i in {0,0.01,...,.8} {\shade[ball color=red]  (1-\i*.574,1+\i*.819) circle (1.5pt);}
    \node[left] at (1-.8*.574,1+.8*.819) {\footnotesize $L$};
\end{tikzpicture}
\end{center}
  \caption{The segment $L$, the planar closed convex cone $C$, the ray $R$ and the point $A$.}\label{cone-ray-point}
\end{figure}
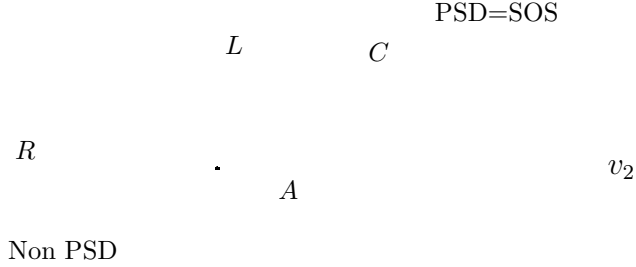

In Section 4, we show that a point $P$ in $S$ is PNS-free if there is a value $M$, such that when $v_0 = M$, $f_0(\x) \equiv f(\x)$ has an SOS decomposition, and $f_0(\bar \x) = 0$ for $\bar \x = (\bar x_1, \bar x_2, \bar x_3, \bar x_4)^\top \in \Re^4$ with $\bar x_1^2+\bar x_4^2 \not = 0$.  We call such a value $M$, such an SOS decomposition of $f_0(\x)$, and such a vector $\bar \x$ the {\bf critical value}, the {\bf critical SOS decomposition} and the {\bf critical minimizer} of $\A$ at $P$, respectively.   Then, we show that the segment $L = \{ (v_2, v_6, v_1, v_3, v_5)^\top = (1, 1, t, t, t)^\top: t \in [-1, 1] \}$ is PNS-free.  We conjecture that this segment is the minimizer set of both $M_0$ and $N_0$.  Then, we show that the $45$-degree planar closed convex cone $C = \{ (v_2, v_6, v_1, v_3, v_5)^\top =(a, b, 0, 0, 0)^\top: a \ge b \ge 1 \}$, the ray $R = \{ (v_2, v_6, v_1, v_3, v_5)^\top = (a, 0, 0, 0, 0)^\top: a \le 0 \}$ and the point  $A = (1, 0, 0, 0, 0)^\top$ are also PNS-free. We illustrate $L$, $C$, $R$ and $A$ in Figure \ref{cone-ray-point}.

In Section 5, numerical tests check various grid points, and find that $M_0 = N_0$ there.   Thus, they are also PNS-free.  Therefore, numerical tests indicate that there are no fourth order four dimensional PNS Hankel tensors with symmetric generating vectors.

Some final remarks are made in Section 6.

\section{Fourth Order Four Dimensional Hankel Tensors}
\hspace{4mm}

We write out (\ref{e1}) explicitly in terms of the coordinates of its generating vector $\vv$:
\begin{eqnarray}
\lefteqn{ f(\x) = v_0x_1^4 + 4v_1x_1^3x_2 + v_2(4x_1^3x_3+6x_1^2x_2^2)
     + v_3(4x_1x_2^3+4x_1^3x_4+12x_1^2x_2x_3) } \nonumber \\
&&{} + v_4(x_2^4+6x_1^2x_3^2+12x_1x_2^2x_3+12x_1^2x_2x_4) + v_5(4x_2^3x_3 + 12x_1x_2x_3^2 + 12x_1x_2^2x_4 + 12x_1^2x_3x_4) \nonumber \\
&&{} + v_6(4x_1x_3^3 + 4x_2^3x_4 + 6x_1^2x_4^2 + 6x_2^2x_3^2 + 24x_1x_2x_3x_4) \label{e6} \\
&&{} + v_7(4x_2x_3^3 + 12x_2^2x_3x_4 + 12x_1x_3^2x_4 + 12x_1x_2x_4^2)
     + v_8(x_3^4+6x_2^2x_4^2+12x_2x_3^2x_4+12x_1x_3x_4^2) \nonumber \\
&&{} + v_9(4x_3^3x_4+4x_1x_4^3+12x_2x_3x_4^2) + v_{10}(4x_2x_4^3+6x_3^2x_4^2) + 4v_{11}x_3x_4^3 + v_{12}x_4^4.  \nonumber
\end{eqnarray}

The following theorem gives some necessary conditions for fourth order four dimensional Hankel tensors being PSD. Particularly, we note that four key elements of its generating vector $v_0,v_4,v_8,v_{12}$ must be nonnegative.

\begin{Theorem} \label{t1}
Suppose that $\A = (a_{i_1i_2i_3i_4})$ is a Hankel tensor generated by
its generating vector $\vv = (v_0, v_1, \ldots, v_{12})^\top \in \Re^{13}$.
If $\A$ is a PSD (or positive definite, or SOS, or strong) Hankel tensor, then we have
\begin{equation} \label{e7}
v_i \ge 0,
\end{equation}
for $i= 0, 4, 8, 12$,
\begin{equation} \label{e8}
v_i + 6v_{i+2} + v_{i+4} \ge 4|v_{i+1} + v_{i+3}|,
\end{equation}
for $i=0, 4, 8$,
\begin{equation} \label{e9}
v_i + 6v_{i+4} + v_{i+8} \ge 4|v_{i+2} + v_{i+6}|,
\end{equation}
for $i=0, 4$, and
\begin{equation} \label{e10}
v_0 + 6v_6 + v_{12} \ge 4|v_3 + v_9|.
\end{equation}
\end{Theorem}

\noindent
{\bf Proof}
Let $\e_k$ be the $k$th column of a $4$-by-$4$ identity matrix, for $k=1,2,3,4$.
Substituting $\x = \e_k$ to (\ref{e6}) for $k=1, 2, 3, 4$, by $f(\e_k) \ge 0$, we have (\ref{e7}) for  $i= 0, 4, 8, 12$.

Substituting $\x = \e_k + \e_{k+1}$ to (\ref{e6}) for $k=1, 2, 3$, by $f(\e_k +\e_{k+1}) \ge 0$, we have
$$v_i +4v_{i+1} + 6v_{i+2} + 4v_{i+3} + v_{i+4} \ge 0,$$
for $i=0, 4, 8$.
Substituting $\x = \e_k - \e_{k+1}$ to (\ref{e6}) for $k=1, 2, 3$, by $f(\e_k -\e_{k+1}) \ge 0$, we have
$$v_i -4v_{i+1} + 6v_{i+2} - 4v_{i+3} + v_{i+4} \ge 0,$$
for $i=0, 4, 8$.   Combining these two inequalities, we have (\ref{e8}) for $i=0, 4, 8$.

Similarly, by $f(\e_k +\e_{k+2}) \ge 0$ and $f(\e_k -\e_{k+2}) \ge 0$ for $k=1, 2$, we have (\ref{e9}) for $i=0, 4$.
By $f(\e_1 +\e_4) \ge 0$ and $f(\e_1 -\e_4) \ge 0$, we have (\ref{e10}).   The theorem is proved.
\ep

Whereafter, we say that a PSD Henkel tensor is SOS if there is a key element of its generating vector $v_0,v_4,v_8,v_{12}$ vanishes. Before we show this, the following lemma is useful.

\begin{Lemma}\label{Lem-oddp}
  If a polynomial in one variable is always nonnegative:
  $$ p(t) = a_0t^{2k+1}+a_1t^{2k}+\cdots+a_{2k+1} \geq 0,
    \qquad \forall\,t\in\Re.  $$
  Then $a_0 = 0.$
\end{Lemma}
\noindent
{\bf Proof}
  If $a_0>0$, we let $t\to-\infty$ and get $p(t)\to-\infty$,
  which contradicts that $p(t)$ is nonnegative.

  If $a_0<0$, we let $t\to+\infty$ and get $p(t)\to-\infty$,
  which also contradicts that $p(t)$ is nonnegative.

  Hence, there must be $a_0=0$.
\ep

\begin{Theorem} \label{t2}
  Suppose the fourth order four dimensional Hankel tensor $\A$ is PSD
  and its generating vector is $\vv$. If $v_0v_{12} = 0$, then $v_j = 0$,
  for $j = 1, \ldots, 11$, and $\A$ is SOS.
\end{Theorem}
\noindent
{\bf Proof}
  Without loss of generality, we assume that $v_0 = 0$.

  To prove $v_1=0$, we take $\x_1=(t,1,0,0)^\top$. Then, the homogeneous polynomial (\ref{e6})
  reduces to
  $$ f(\x_1) = 4v_1t^3 + 6v_2t^2 + 4v_3t +v_4. $$
  From Lemma \ref{Lem-oddp}, we have $v_1=0$ since $f(\x_1)$ is nonnegative.
  Similarly, we can prove $v_2=v_3=0$ if we take $\x_2=(t,0,1,0)^\top$ and $\x_3=(t,0,0,1)^\top$ respectively.

  From Theorem \ref{t1}, we know $v_4 \geq 0$.
  When we take $\x_4=(t^2,t,-\frac{1}{\sqrt{6}},0)^\top$,
  the homogeneous polynomial (\ref{e6}) reduces to
  $$ f(\x_4) = -(2\sqrt{6}-2)v_4t^4 + \mathcal{O}(t^3). $$
  Let $t\to\infty$.
  Since $f(\x_4)$ is always nonnegative, we have $v_4 \leq 0$.
  Hence, there must be $v_4=0$.

  If we take $\x_5=(t^3,0,t,1)^\top$, the homogeneous polynomial (\ref{e6}) is
  $$ f(\x_5) = 12v_5t^7 + \mathcal{O}(t^6).$$
  From Lemma \ref{Lem-oddp}, we have $v_5=0$ since $f(\x_5)$ is nonnegative.

  We take $\x_6=(t,0,1,0)^\top$. Then,
  the homogeneous polynomial (\ref{e6}) is
  $$ f(\x_6) = 4v_6t + v_8. $$
  From Lemma \ref{Lem-oddp}, we have $v_6=0$ since $f(\x_6)$ is nonnegative.
  Similarly, we can prove $v_7=0$ when we take $\x_7=(0,t,1,0)^\top$.

  We take $\x_8=(t^4,0,t,1)^\top$. Then we have
  $$ f(\x_8) = 12v_8t^5 + \mathcal{O}(t^4).$$
  From Lemma \ref{Lem-oddp}, we have $v_8=0$ since the polynomial $f(\x_8)$ is nonnegative.

  We could prove $v_9=0$, $v_{10}=0$ and $v_{11}=0$ if we takes
  $\x_9=(t,0,0,1)^\top$, $\x_{10}=(0,t,0,1)^\top$ and $\x_{11}=(0,0,t,1)^\top$, respectively.

  Finally, since $v_0=v_1=\cdots=v_{11}=0$, we have
  $$f(\x) = v_{12}x_4^4.$$
  By Theorem \ref{t1}, we get $v_{12}\geq 0$.
  Hence, the Hankel tensor $\A$ is obviously SOS.
\ep

\begin{Theorem} \label{t3}
  Suppose the fourth order four dimensional Hankel tensor $\A$ is PSD
  and its generating vector is $\vv$. If $v_4v_8 = 0$, then $v_j = 0$
  for $j = 1, 2, \ldots, 11$, and $\A$ is SOS.
\end{Theorem}
\noindent
{\bf Proof}
  By symmetry, we only need to prove this theorem under the condition $v_4 = 0$.

  If we take $\x_1=(1,t,0,0)^\top$,
  the homogeneous polynomial (\ref{e6}) reduces to
  $$ f(\x_1) = 4v_3t^3 + 6v_2t^2 + 4v_1t +v_0. $$
  From Lemma \ref{Lem-oddp}, we have $v_3=0$ since $f(\x_1)$ is nonnegative.
  Similarly, we can prove $v_5=v_6=0$ if we take $\x_2=(0,t,1,0)^\top$ and $\x_3=(0,t,0,1)^\top$ respectively.

  To prove $v_7=0$, we take $\x_4=(0,t^2,t,1)^\top$. Then, the homogeneous polynomial (\ref{e6})
  reduces to
  $$ f(\x_4) = 16v_7t^5 + \mathcal{O}(t^4). $$
  From Lemma \ref{Lem-oddp}, we have $v_7=0$ since $f(\x_4)$ is nonnegative.

  From Theorem \ref{t1}, we know $v_8 \geq 0$.
  When we take $\x_5=(0,-t^2,t,1)^\top$,
  the homogeneous polynomial (\ref{e6}) reduces to
  $$ f(\x_5) = -5v_8t^4 + \mathcal{O}(t^3). $$
  Let $t\to\infty$.
  Since $f(\x_5)$ is always nonnegative, we have $v_8 \leq 0$.
  Hence, there must be $v_8=0$.

  If we take $\x_6=(0,0,t,1)^\top$, the homogeneous polynomial (\ref{e6}) is
  $$ f(\x_6) = 4v_9t^3 + \mathcal{O}(t^2).$$
  From Lemma \ref{Lem-oddp}, we have $v_9=0$ since $f(\x_6)$ is nonnegative.
  Similarly, we could prove $v_{10}=0$ and $v_{11}=0$ if we takes
  $\x_7=(0,t,0,1)^\top$ and $\x_8=(0,0,t,1)^\top$, respectively.

  The prove of $v_1=0$ and $v_2=0$ could be similarly obtained if we take
  $\x_9=(1,t,0,0)^\top$ and $\x_{10}=(1,0,t,0)^\top$ respectively.

  Finally, since $v_j = 0$ for $j = 1,\ldots,11$, we have
  $$f(\x) = v_0x_1^4 + v_{12}x_4^4.$$
  By Theorem \ref{t1}, we get $v_0\geq 0$ and $v_{12}\geq 0$.
  Hence, the Hankel tensor $\A$ is obviously SOS.
\ep



\section{Symmetric Generating Vectors}
\hspace{4mm}

Now, we make assumptions (\ref{e2}) and (\ref{e3}).
At the beginning, we consider a mini problem which is the Hankel polynomial with $x_1=x_4=0$.
This problem helps us to analyze the effective domain of two important surfaces $M_0$ and $N_0$.



\subsection{Function $\eta$}
\hspace{4mm}

We consider a two variable quartic polynomial
$$g(y_1, y_2) = \alpha y_1^4 + 4\beta y_1^3y_2 + 6\gamma y_1^2y_2^2 + 4\beta y_1y_2^3 + \alpha y_2^4.$$
Its PSD property is completely characterized by the following theorem.

\begin{Theorem} \label{p-1}
  The quartic polynomial $g(y_1,y_2)$ is PSD if and only if
  \begin{equation*}
    \alpha\ge \eta(\beta,\gamma) := \left\{\begin{aligned}
      4|\beta | - 3\gamma                         ~&~ \text{ if } \gamma\le |\beta|,\\
      \frac{3\gamma-\sqrt{9\gamma^2-8\beta^2}}{2} ~&~ \text{ if } \gamma> |\beta|. \\
    \end{aligned}\right.
  \end{equation*}
\end{Theorem}
\noindent
{\bf Proof}
  First, if $g(y_1,y_2)$ is PSD, from $g(1, -1) \ge 0$ and $g(1, 1) \ge 0$,
  we have $\alpha \ge 4|\beta | - 3\gamma$.  Thus, in any case,
  $\eta(\beta, \gamma) \ge 4|\beta | - 3\gamma$.

  Second, suppose that $\alpha \ge 4|\beta | - 3\gamma$.   If $\gamma \le 0$, we get
  $$g(y_1, y_2) = (\alpha - 4|\beta | + 3\gamma)(y_1^4+ y_2^4) + 4|\beta |(y_1+y_2)^2(y_1^2-y_1y_2+y_2^2) -3\gamma(y_1^2-y_2^2)^2 \ge 0.$$
  If $0<\gamma \le |\beta|$, we rewrite $g(y_1,y_2)$ as follows
  $$g(y_1, y_2) = (\alpha - 4|\beta | + 3\gamma)(y_1^4+ y_2^4) + (y_1+y_2)^2\left[(4|\beta| - 3\gamma)(y_1^2+y_2^2)-(4|\beta|-6\gamma)y_1y_2\right].$$
  Since $(4|\beta|-6\gamma)^2-4(4|\beta|-3\gamma)^2=-48|\beta|(|\beta|-\gamma)\le0$, it yields that $g(y_1,y_2)\ge 0$.

  Finally, we consider the case $\gamma > |\beta|$.
  Let $\bar{\alpha}=\frac{3\gamma-\sqrt{9\gamma^2-8\beta^2}}{2} >0$. Then, we have
  \begin{equation*}
  g(y_1,y_2) = (\alpha-\bar{\alpha})(y_1^4+y_2^4)+
      \bar{\alpha}\left(y_1^2+\frac{2\beta}{\bar{\alpha}}y_1y_2+y_2^2\right)^2.
  \end{equation*}
  Obviously, if $\alpha\geq\bar{\alpha}$, $g(y_1,y_2)$ is SOS and PSD.

  Next, we show that $y_1^2+\frac{2\beta}{\bar{\alpha}}y_1y_2+y_2^2 = 0$ has nonzero real roots.
  For the convenience, we denote $t=\frac{y_1}{y_2}$ and prove
  that $t^2+\frac{2\beta}{\bar{\alpha}}t+1 = 0$ has real roots.
  It is easy to see that $t=0$ is not its root.
  Since $\gamma>|\beta|$, we have
  \begin{equation*}
    \frac{|\beta|}{\bar{\alpha}} = \frac{2|\beta|}{3\gamma-\sqrt{9\gamma^2-8\beta^2}}
     = \frac{2|\beta|(3\gamma+\sqrt{9\gamma^2-8\beta^2})}{8\beta^2}
     \geq \frac{8|\beta|\gamma}{8\beta^2}
     \geq 1.
  \end{equation*}
  Hence, $|\beta| \geq \bar{\alpha}$.
  The discriminant of the quadratic in $t$ is
  $$ \left(\frac{2\beta}{\bar{\alpha}}\right)^2-4 = 4\frac{\beta^2-\bar{\alpha}^2}{\bar{\alpha}^2} \geq 0. $$
  Therefore, there are nonzero $(y_1,y_2)$ such that $g(y_1,y_2)=(\alpha-\bar{\alpha})(y_1^4+y_2^4)$. Obviously, if $g(y_1,y_2)$ is PSD, we have $\alpha\ge\bar{\alpha}$.
  Thus, we say $\eta(\beta,\gamma)=\bar{\alpha}$ if $\gamma>|\beta|$.
\ep

Then we have another necessary condition for a fourth order four dimensional Hankel tensor $\A$ to be PSD under assumptions (\ref{e2}) and (\ref{e3}).

\begin{Corollary} \label{p0}
Under assumptions (\ref{e2}) and (\ref{e3}), if $\A$ is PSD, then $\eta(v_5, v_6) \le 1$.
\end{Corollary}

\noindent
{\bf Proof}  Let $x_1 = x_4 = 0$, $x_2 = y_1$ and $x_3 = y_2$.  By Theorem \ref{p-1}, we have the conclusion.  \ep

\subsection{Surfaces $M_0$ and $N_0$}
\hspace{4mm}

We now introduce the key idea of this paper, to establish two surface $M_0$ and $N_0$, in the following theorem.

\begin{Theorem}  \label{t5}
Suppose that assumptions (\ref{e2}) and (\ref{e3}) hold.
Then, there are two functions $M_0(v_2, v_6, v_1, v_3, v_5) \geq N_0(v_2, v_6, v_1, v_3, v_5) > 0$ defined for
\begin{equation}  \label{e12}
\eta(v_5, v_6) < 1,
\end{equation}
such that $\A$ is SOS if and only if $v_0 \ge M_0(v_2, v_6, v_1, v_3, v_5)$,
and $\A$ is PSD if and only if $v_0 \ge N_0(v_2, v_6, v_1, v_3, v_5)$.  If for all $v_5$ and $v_6$ satisfying (\ref{e12}), we have
$M_0(v_2, v_6, v_1, v_3, v_5) = N_0(v_2, v_6, v_1, v_3, v_5)$, then there are no fourth order four dimensional PNS Hankel tensors under assumption (\ref{e2}).
\end{Theorem}

\noindent
{\bf Proof}   Using assumptions (\ref{e2}) and (\ref{e3}), we rewrite (\ref{e6}) as
$$f(\x) = v_0(x_1^4+x_4^4) + \bar v_4(x_2^4 + x_3^4) + f_1(\x) + f_2(\x),$$
where
$$f_1(\x) = \eta(v_5, v_6)(x_2^4 + x_3^4) + 4v_5(x_2^3x_2+x_2x_3^3) + 6v_6x_2^2x_3^2$$
and
$$\bar v_4 = 1 - \eta(v_5, v_6).$$
Then $\bar v_4 > 0$ by (\ref{e12}).    By Theorem \ref{p-1}, $f_1(\x)$ is PSD.   Since $f_1(\x)$ has only two variables, it is also SOS by Hilbert \cite{Hi, Re}.

We now consider terms in $f_2(\x)$.  Each monomial in $f_2(\x)$ has at least one factor as a power of $x_1$ or $x_4$.  We may order the monomials of $f_2(\x)$.   For example, consider $12v_5x_1x_2x_3^2$.   Assume that it is ordered as the $k$th monomial of $f_2(\x)$.
Then by the arithmetic-geometric inequality, we may see that
$$-12v_5x_1x_2x_3^2 \le 3|v_5|\left({1 \over \epsilon_k^3}x_1^4 + \epsilon_kx_2^4 + 2\epsilon_kx_3^4\right),$$
where $\epsilon_k$ is a small positive number.  We may let $\epsilon_k$ be small enough such that the sum of the coefficients for $x_2^4$ on the right hand side of the above inequality for all possible $k$ is less than $\bar v_4$.   By symmetry, the sum of the coefficients for $x_3^4$ on the right hand side of the above inequality for all possible $k$ is less than $\bar v_4$.   We see that
$$12v_5x_1x_2x_3^2 + 3|v_5|\left({1 \over \epsilon_k^3}x_1^4 + \epsilon_kx_2^4 + 2\epsilon_kx_3^4\right)$$
is a PSD diagonal minus tail form.   By \cite{FK}, it is SOS.   Thus, as long as $v_0$ is big enough, when (\ref{e12}) is satisfied, $f(\x)$ is SOS.  From this, we see that $M_0$ and $N_0$ exist, such that they are defined as long as (\ref{e12}) is satisfied, $M_0 \ge N_0$, $\A$ is SOS if and only if $v_0 \ge M_0$, and
$\A$ is PSD if and only if $v_0 \ge N_0$.

By Theorem \ref{p-1}, we now only need to consider the case that $\eta(v_5, v_6) = 1$.   Suppose that for all $v_5$ and $v_6$ satisfying (\ref{e12}), we have
$M_0(v_2, v_6, v_1, v_3, v_5) = N_0(v_2, v_6, v_1, v_3, v_5)$.  Since the sets for PSD Hankel tensors and SOS Hankel tensors are closed \cite{LQX}, this implies that for all $v_5$ and $v_6$ satisfying $\eta(v_5, v_6) = 1$, we also have
$M_0(v_2, v_6, v_1, v_3, v_5) = N_0(v_2, v_6, v_1, v_3, v_5)$, as long as $N_0$ is defined there.   Thus, in this case, by Theorem \ref{t3}, there are no fourth order four dimensional PNS Hankel tensors under assumption (\ref{e2}).
\ep

For the variables of $M_0$ and $N_0$, we put $v_2$ and $v_6$ before $v_1$, $v_3$ and $v_5$, as $v_2, v_6$ play a more important role in the PSD and SOS properties of $\A$, comparing with $v_1, v_3$ and $v_5$.

We now regard $P = (v_2, v_6, v_1, v_3, v_5)^\top$ as a point in $\Re^5$.   If $M_0(P) = N_0(P)$, $P$ is called a {\bf PNS-free point}.  We call
$$S = \{ (v_2, v_6, v_1, v_3, v_5)^\top \in \Re^5 : \eta(v_5, v_6) < 1 \}$$
the {\bf effective domain}.   Theorem \ref{t5} says that if all the points in the effective domain are PNS-free, then there are no fourth order four dimensional PNS Hankel tensors with symmetric generating vectors.    In the next sections, we will study more on PNS-free points.

\section{Theoretical Proofs of Some PNS-Free Regions}
\hspace{4mm}

\subsection{Critical SOS Decomposition}
\hspace{4mm}

For the convenience, we present formally three ingredients used in theoretical proofs of this section.
If a point belongs to the effective domain and enjoys these ingredients, it is PNS-free.

\begin{Definition}
Suppose that assumptions (\ref{e2}) and (\ref{e3}) hold and $P = (v_2, v_6, v_1, v_3, v_5)^\top \in S$.
Suppose that there is a number $M$ such that $\A$ is SOS if $v_0 = M$, and
a point $\bar \x = (\bar x_1, \bar x_2, \bar x_3, \bar x_4)^\top \in \Re^4$ such that
$\bar x_1^2 + \bar x_4^2 > 0$ and $f_0(\bar \x) = 0$, where $f_0(\x) \equiv f(\x)$ with $v_0 = M$.
Then we call $M$ the {\bf critical value} of $\A$ at $P$, the SOS decomposition $f_0(\x)$ the {\bf critical SOS decomposition} of $\A$ at $P$, and $\bar \x$ the {\bf critical minimizer} of $\A$ at $P$.
\end{Definition}

\begin{Theorem} \label{t6}
Let $P \in S$.   Then $P$ is PNS-free if $\A$ has a critical value $M$, a critical SOS decomposition $f_0(\x)$ and a critical minimizer $\bar{\x}$ at $P$.
\end{Theorem}

\noindent
{\bf Proof}  Suppose that $\A$ has a critical value $M$, a critical SOS decomposition $f_0(\x)$ and a critical minimizer $\bar{\x}$ at $P$.   Then we have $M \ge M_0(P)$ by the definition of $M_0$.   If $v_0 < M$, then
$$f(\bar \x) = (v_0-M)(\bar x_1^4+\bar x_4^4) + f_0(\bar \x) < 0.$$
This implies that $N_0(P) \ge M$ by the definition of $N_0$.  But $N_0(P) \le M_0(P)$.  Thus, $M_0(P) = N_0(P) = M$, i.e., $P$ is PNS-free.
\ep

We believe that all the effective domain $S$ is PNS-free.  In the next four subsections, we theoretically prove that some regions of $S$ are PNS-free.

\subsection{A PNS-Free Segment}
\hspace{4mm}

Professor Man Kam Kwong pointed out that $N_0(1, 1, 0, 0, 0) = 1$, $N_0(2, 1, 0, 0, 0) = 8$ and $N_0(4, 0, 0, 0, 0) = 441$, are integers. See also Table 1 in Section 6.   He suggested us to considered these three points more carefully.     Stimulated by Prof. Kwong's comments, we derive the results of Subsections 4.2 and 4.3.

We have the following theorem.

\begin{Theorem} \label{t21}
Suppose that $P = (v_2, v_6, v_1, v_3, v_5)^\top = (1, 1, t, t, t)^\top$, where $t \in [-1, 1]$.  Then, $P$ is PNS-free, with the critical value $1$ and the critical minimizer $(1, 0, -1, 0)^\top$.
\end{Theorem}

\noindent
{\bf Proof}
For $P = (v_2, v_6, v_1, v_3, v_5)^\top = (1, 1, t, t, t)^\top$, where $t \in [-1, 1]$, and $M = 1$, we have
$$f_0(\x) = {1+t \over 2}(x_1+x_2+x_3+x_4)^4 +{1-t \over 2}(x_1-x_2+x_3-x_4)^4$$
is SOS, and
$$f_0(1, 0, -1, 0) = 0.$$
Hence, $P$ is PNS-free.
\ep

By numerical experiments, we have the following conjecture.

\begin{Conjecture}
The segment $L = \{ (v_2, v_6, v_1, v_3, v_5)^\top = (1, 1, t, t, t)^\top: t \in [-1, 1] \}$, is the minimizer set of both $M_0$ and $N_0$.
\end{Conjecture}

\subsection{A PNS-free Planar Cone}
\hspace{4mm}

\begin{Theorem} \label{t-bgeOne}
Suppose that $P = (v_2, v_6, v_1, v_3, v_5)^\top = (v_2, v_6, 0, 0, 0)^\top$ with $v_2 \ge v_6 \ge 1$.
Then, $P$ is PNS-free.

If we parameterize $v_6=b$ and $v_2=(\theta+3b-1)(\theta^2+(3b-2)\theta-3b+4)$.
Then, the critical value at $P$ is
\begin{equation*}
    M=(\theta+3b-1)^2(3\theta^2+(10b-6)\theta+3b^2-10b+9)
\end{equation*}
and the critical minimizer is $\bar{\x}=(1, 0, -(\theta+3b-1), 0)^\top$.
\end{Theorem}

\noindent
{\bf Proof}   Note that for $v_2 \ge v_6 \ge 1$, we may let $v_6 = b$ and
$v_2 = (\theta+3b-1)(\theta^2+(3b-2)\theta-3b+4)$, where the parameter
$$\theta \ge \bar{\theta}=(b-1)^{1 \over 3}(b+1)^{2 \over 3} + (b-1)^{2 \over 3}(b+1)^{1 \over 3} -2b+1.$$
In fact, $\bar{\theta}$ is the largest real root of the cubic equation $v_2-v_6=0$.

With the critical value as $M=(\theta+3b-1)^2(3\theta^2+(10b-6)\theta+3b^2-10b+9)$, the critical SOS decomposition at $P$ is as follows
  \begin{eqnarray*}
    f_0(\x) &=& \frac{1}{v_0}(v_0x_1^2+2v_2x_1x_3+\alpha_1x_3^2)^2
              + \frac{1}{v_0}(v_0x_4^2+2v_2x_2x_4+\alpha_1x_2^2)^2 \\
      &&{} + \alpha_2((\theta+3b-1)x_1x_3+x_3^2)^2 + \alpha_2((\theta+3b-1)x_2x_4+x_2^2)^2 \\
      &&{} + \frac{6}{b}(x_1x_2+x_3x_4+bx_2x_3+bx_1x_4)^2 + \frac{6(b^2-1)}{b}(x_1x_2 + x_3x_4)^2 \\
      &&{} + 6(v_2-b)[x_1^2x_2^2+x_3^2x_4^2],
  \end{eqnarray*}
  where the involved parameters are as follows:
  \begin{eqnarray*}
    \alpha_1 &=& -(\theta^2+(4b-2)\theta+3b^2-4b+1), \\
    \alpha_2 &=& \frac{2(\theta^2+(4b-2)\theta+b^2-4b+4)}{3\theta^2+(10b-6)\theta+3b^2-10b+9}.
  \end{eqnarray*}

  Since $f_0(1,0,-(\theta+3b-1),0)=0$,
  the corresponding critical minimizer is $\bar{\x}=(1,0,-(\theta+3b-1),0)^{\top}$.
  Hence, $P =(v_2,v_6, 0, 0, 0)^\top$
  with $v_2 \ge v_6 \ge 1$ is PNS-free.
\ep

The cone $C = \{ (v_2, v_6, v_1, v_3, v_5)^\top =(a, b, 0, 0, 0)^\top: a \ge b \ge 1 \}$ is a $45$-degree planar closed convex cone.  Its end point is just the mid point of the segment $L = \{ (v_2, v_6, v_1, v_3, v_5)^\top = (1, 1, t, t, t)^\top: t \in [-1, 1] \}$, discussed in the last subsection.

\subsection{A PNS-Free Ray}
\hspace{4mm}

In this subsection, we show that the ray $R = \{ (v_2, v_6, v_1, v_3, v_5)^\top = (a, 0, 0, 0, 0)^\top: a \le 0 \}$ is PNS-free. Let $a=-\rho$, where $\rho\ge 0$ is a constant.
We report that, at a point $P=(-\rho,0,0,0,0)^{\top}$, $\A$ has the critical value
\begin{equation*}
    M=3\sqrt[3]{\theta_1+32\sqrt{\theta_2}}+\frac{\theta_3}{3\sqrt[3]{\theta_1+32\sqrt{\theta_2}}}
      +6\rho^2+138\rho+609,
\end{equation*}
where
\begin{eqnarray*}
  \theta_1 &:=& -\rho^6+272\rho^5+12608\rho^4+204032\rho^3+1558528\rho^2+5750784\rho+8290304, \\
  \theta_2 &:=& -(\rho+6)^2(\rho+4)^3(\rho^2+4\rho-16)^3, \\
  \theta_3 &:=& 9(\rho+8)(\rho^3+152\rho^2+1728\rho+5120).
\end{eqnarray*}
The function $f_0(\x)$ enjoys a critical SOS decomposition:
\begin{equation*}
  f_0(\x) = \sum_{k=1}^5 q_k^2(\x),
\end{equation*}
where
\begin{eqnarray*}
  q_1(\x) &=& x_3^2+6x_2x_4+\alpha_1x_1^2+\alpha_2x_4^2, \\
  q_2(\x) &=& x_2^2+6x_1x_3+\alpha_2x_1^2+\alpha_1x_4^2, \\
  q_3(\x) &=& \alpha_3x_2x_4+\alpha_4x_1^2+\alpha_5x_4^2, \\
  q_4(\x) &=& \alpha_3x_1x_3+\alpha_5x_1^2+\alpha_4x_4^2, \\
  q_5(\x) &=& \alpha_6x_1^2-\alpha_6x_4^2.
\end{eqnarray*}
The involved parameters are listed as follows:
\begin{eqnarray*}
  \alpha_1 &=& -\frac{(\rho+23)M_1(-\rho)-9\rho^3-21\rho^2+105\rho+9}{M_1(-\rho)+3\rho^2+6\rho-33}, \\
  \alpha_2 &=& -3\rho, \\
  \alpha_3 &=& \sqrt{-30-2\alpha_{15}}, \\
  \alpha_4 &=& \frac{6(1-\alpha_{15})}{\alpha_{33}}, \\
  \alpha_5 &=& \frac{16\rho}{\alpha_{33}}, \\
  \alpha_6 &=& \sqrt{-6\rho\alpha_{15} -\frac{192\rho(\alpha_{15}-1)}{\alpha_{33}^2}}.
\end{eqnarray*}

\begin{Theorem} \label{t10}
Suppose that assumptions (\ref{e2}) and (\ref{e3}) hold.
Then, for any constant $\rho \ge 0$, $P=(-\rho, 0, 0, 0, 0)^\top$ is PNS-free.
\end{Theorem}

\noindent
{\bf Proof}
We only need to prove that there is a critical minimizer. Let $$ \bar{\x} = (\alpha_{33},\alpha_{35}+\alpha_{36}, -\alpha_{35}-\alpha_{36},-\alpha_{33})^{\top}. $$
Then, we get $q_3(\bar{\x})=q_4(\bar{\x})=q_5(\bar{\x})=0$ immediately.
Moreover, we have
\begin{equation*}
  q_1(\bar{\x}) = q_2(\bar{\x})= (\alpha_{35}+\alpha_{36})^2 -6(\alpha_{35}+\alpha_{36})\alpha_{33}+\alpha_{15}\alpha_{33}^2-3\rho\alpha_{33}^2 \\
    = 0.
\end{equation*}
We check the validation of the last equality by a mathematical software Maple. Hence, $f_0(\bar{\x})=0$ and $\bar{\x}$ is a critical minimizer at $P$. Hence, we get the conclusion by Theorem \ref{t6}.
\ep

\subsection{A PNS-Free Point }
\hspace{4mm}

We now show that the point $A=(1, 0, 0, 0, 0)^\top$ is PNS-free.
In fact, the critical value at $A$ is
$$ M=477+3\sqrt[3]{3906351+9120\sqrt{57}}+\frac{74403}{\sqrt[3]{3906351+9120\sqrt{57}}}. $$
The critical SOS decomposition of $f_0(\x)$ is as follows
\begin{equation*}
    f_0(\x)=\sum_{k=1}^7 q_k(\x)^2,
\end{equation*}
where
\begin{eqnarray*}
  q_1(\x) &=& x_3^2+6x_2x_4-21x_1^2+\alpha_1x_4^2, \\
  q_2(\x) &=& x_2^2+6x_1x_3-21x_4^2+\alpha_1x_1^2, \\
  q_3(\x) &=& {2\sqrt{3}}x_2x_4+\alpha_2x_1^2+\alpha_3x_4^2, \\
  q_4(\x) &=& {2\sqrt{3}}x_1x_3+\alpha_2x_4^2+\alpha_3x_1^2, \\
  q_5(\x) &=& \alpha_4x_1^2-\alpha_4x_4^2, \\
  q_6(\x) &=& \beta_1x_1x_2+\beta_2x_1x_4, \\
  q_7(\x) &=& \beta_1x_3x_4+\beta_2x_1x_4.
\end{eqnarray*}
Some involved parameters are listed as follows:
\begin{eqnarray*}
  \beta_1 &=& \frac{\sqrt{-6(M_2-36)(3M_2-4336)}}{\sqrt{M_2^2-1302M_2+25056}}, \\
  \beta_2 &=& \frac{\beta_1(3\beta_1^2+116)}{\beta_1^2+12}, \\
  \alpha_1 &=& 3-\frac{1}{2}\beta_1^2, \\
  \alpha_2 &=& 22\sqrt{3}-\frac{\sqrt{3}}{6}\beta_1\beta_2, \\
  \alpha_3 &=& -\frac{8\sqrt{3}}{3}+\frac{\sqrt{3}}{2}\beta_1^2, \\
  \alpha_4 &=& \sqrt{-42\alpha_1+2\alpha_2\alpha_3+\beta_2^2}.
\end{eqnarray*}

\begin{Theorem}
Suppose that assumptions (\ref{e2}) and (\ref{e3}) hold.
Then, $A=(1, 0, 0, 0, 0)^\top$ is PNS-free.
\end{Theorem}

\noindent
{\bf Proof}
Using the mathematical software Maple, we calculate
\begin{equation*}
    f(\x)-\sum_{k=1}^7 q_k^2(\x) = \frac{-\beta_1^6-120\beta_1^4+(4v_0-4944)\beta_1^2+48v_0-69376}{4(\beta_1^2+12)}
        (x_1^4+x_4^4).
\end{equation*}
Substituting the value of $v_0=M$ and $\beta_1$, we get $f_0(\x)-\sum_{k=1}^7 q_k^2(\x) = 0$.

Let $\bar{\x}=(\beta_1,\beta_2,-\beta_2,-\beta_1)^{\top}$.
Obviously, we obtain $q_5(\bar{\x})=q_6(\bar{\x})=q_7(\bar{\x})=0$.
We find that $q_3(\bar{\x})$ and $q_4(\bar{\x})$ vanishes
if we rewrite all the parameters using $\beta_1$.
Using the value of each parameter, we find that $q_1(\bar{\x})=q_2(\bar{\x})=0$.
Since $\bar{x}_1=\beta_1 \approx 1.73$, $\bar{\x}$ is the critical minimizer.
Therefore, this theorem is valid according to  Theorem \ref{t6}.
\ep

\section{Numerical Experiments} \label{NumerE}
\hspace{4mm}

\begin{table}[bt]
  \centering
  \begin{tabular}{c|cccc cccc}
    \hline\hline
    $v_2~\backslash~v_6$ & $-.2$ & $-.1$ & $0$ & $.5$ & $1$ & $1.5$ & $2$ & $4$ \\
    \hline
-4.0   &   3.54e4 &   8.74e3 &   3.76e3 &   4.78e2 &   3.12e2 &   3.92e2 &   6.23e2 &   6.37e3 \\
-2.0   &   2.98e4 &   6.77e3 &   2.73e3 &   2.75e2 &   1.25e2 &   1.70e2 &   3.57e2 &   6.11e3 \\
-1.0   &   2.72e4 &   5.85e3 &   2.26e3 &   1.91e2 &   6.15e1 &   9.26e1 &   2.73e2 &   6.06e3 \\
-0.5   &   2.59e4 &   5.42e3 &   2.04e3 &   1.53e2 &   3.78e1 &   6.41e1 &   2.48e2 &   6.06e3 \\
 0.0   &   2.46e4 &   4.99e3 &   1.82e3 &   1.20e2 &   1.96e1 &   4.50e1 &   2.39e2 &   6.07e3 \\
 0.5   &   2.34e4 &   4.57e3 &   1.62e3 &   8.90e1 &    7.058 &   4.18e1 &   2.45e2 &   6.09e3 \\
 1.0   &   2.21e4 &   4.17e3 &   1.42e3 &   6.21e1 &    1.000 &   4.93e1 &   2.56e2 &   6.11e3 \\
 1.5   &   2.09e4 &   3.78e3 &   1.23e3 &   3.90e1 &    4.191 &   5.69e1 &   2.67e2 &   6.14e3 \\
 2.0   &   1.98e4 &   3.41e3 &   1.06e3 &   2.02e1 &   8.00e0 &   6.46e1 &   2.78e2 &   6.16e3 \\
 3.0   &   1.75e4 &   2.70e3 &   7.28e2 &   7.16e0 &   1.66e1 &   8.01e1 &   3.01e2 &   6.21e3 \\
 4.0   &   1.53e4 &   2.04e3 &   4.41e2 &   1.23e1 &   2.60e1 &   9.60e1 &   3.23e2 &   6.25e3 \\
    \hline\hline
  \end{tabular}
  \caption{The values of $M_0(v_2,v_6,0,0,0)=N_0(v_2,v_6,0,0,0)$ on some grid points.}\label{Table-1}
\end{table}

\begin{figure}[tb]
  \centering
  \includegraphics[width=.7\textwidth]{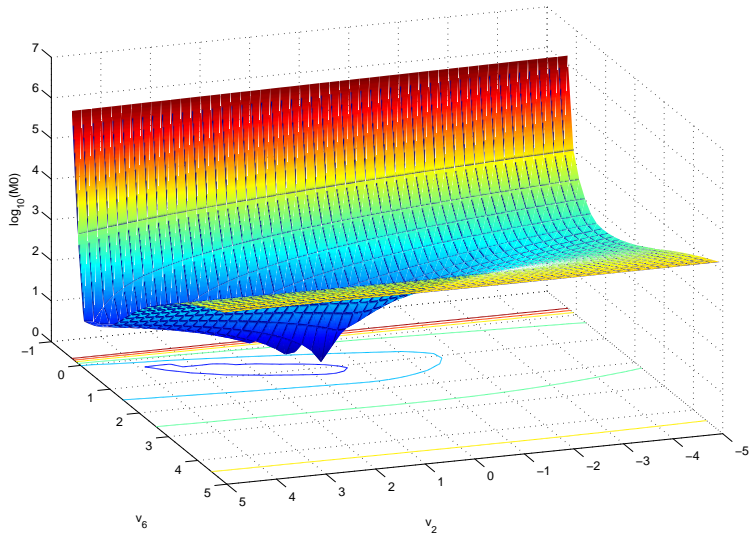}\\
  \caption{The contour profile of $M_0(v_2,v_6,0,0,0)=N_0(v_2,v_6,0,0,0)$.}\label{fig01}
\end{figure}

\begin{figure}[tb]
  \centering
  \includegraphics[width=.7\textwidth]{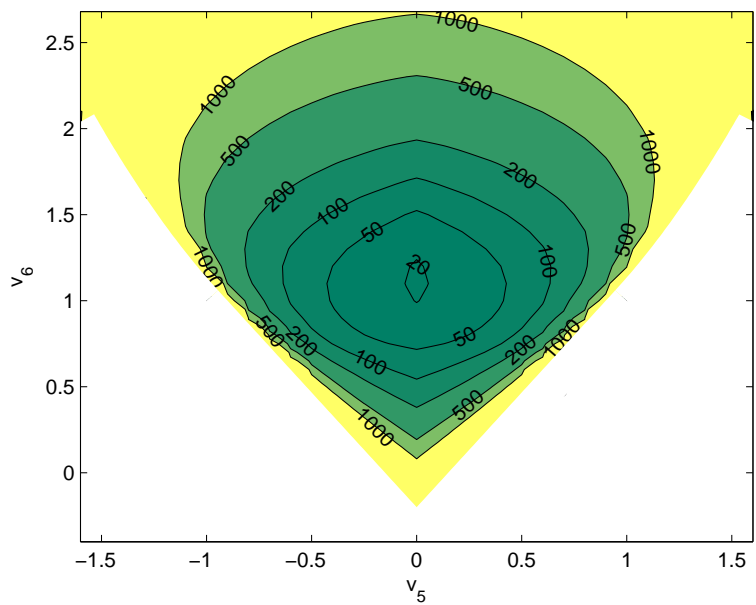}\\
  \caption{The contour profile of $M_0(0,v_6,0,0,v_5)$.}\label{fig02}
\end{figure}

\begin{figure}[tb]
\begin{center}
\begin{tabular}{cc}
  \includegraphics[width=.48\textwidth]{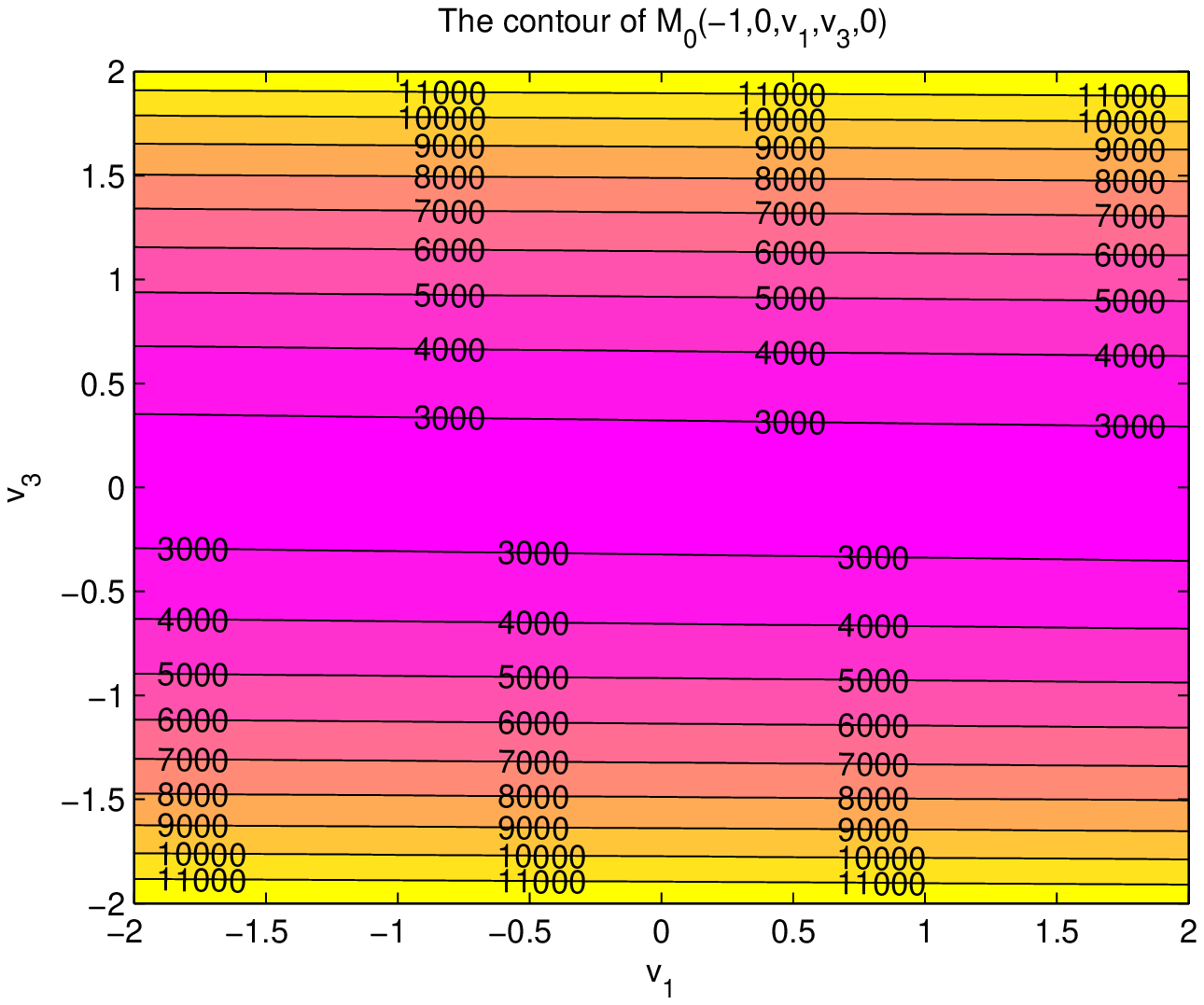} &
  \includegraphics[width=.48\textwidth]{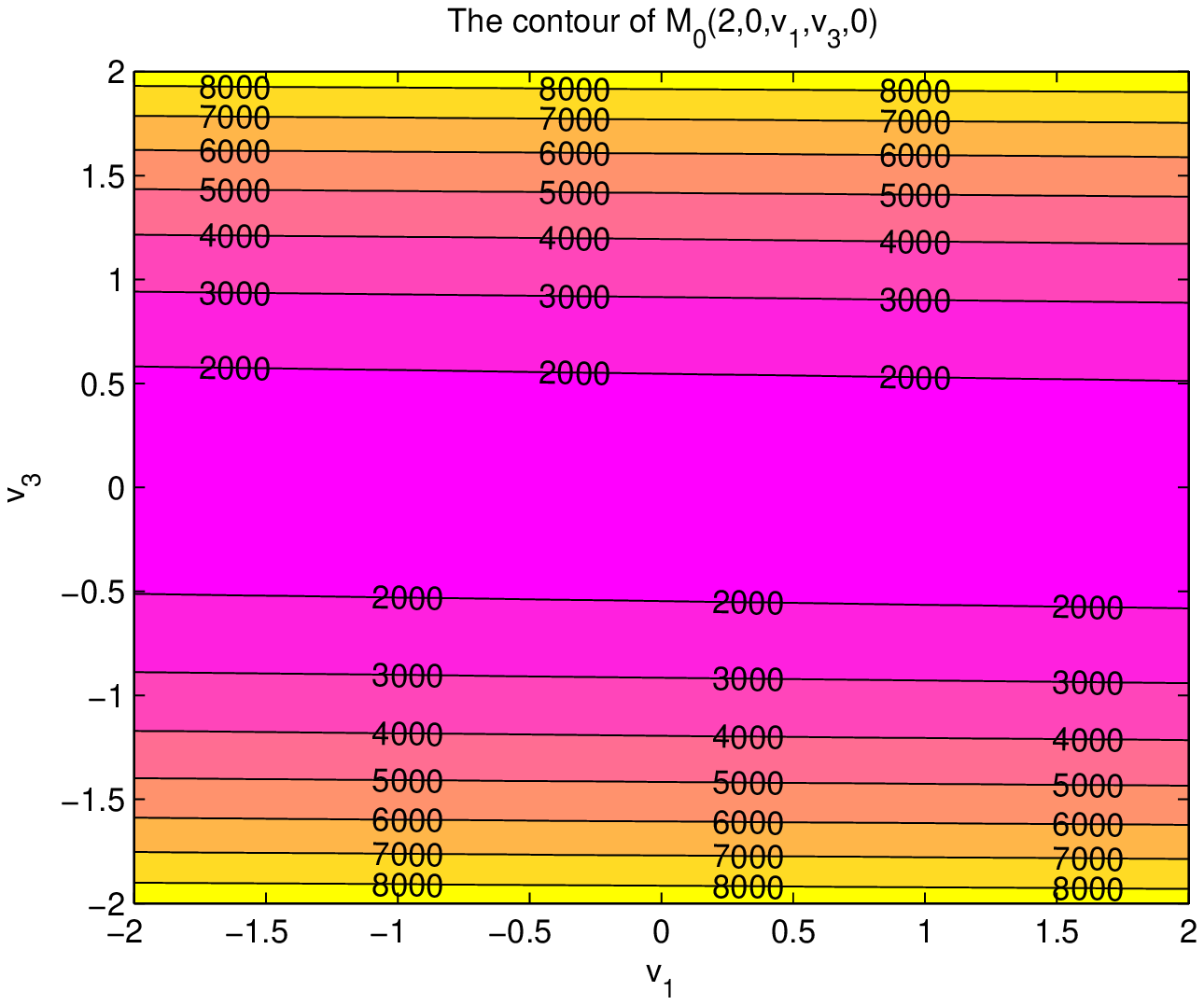} \\
  \includegraphics[width=.48\textwidth]{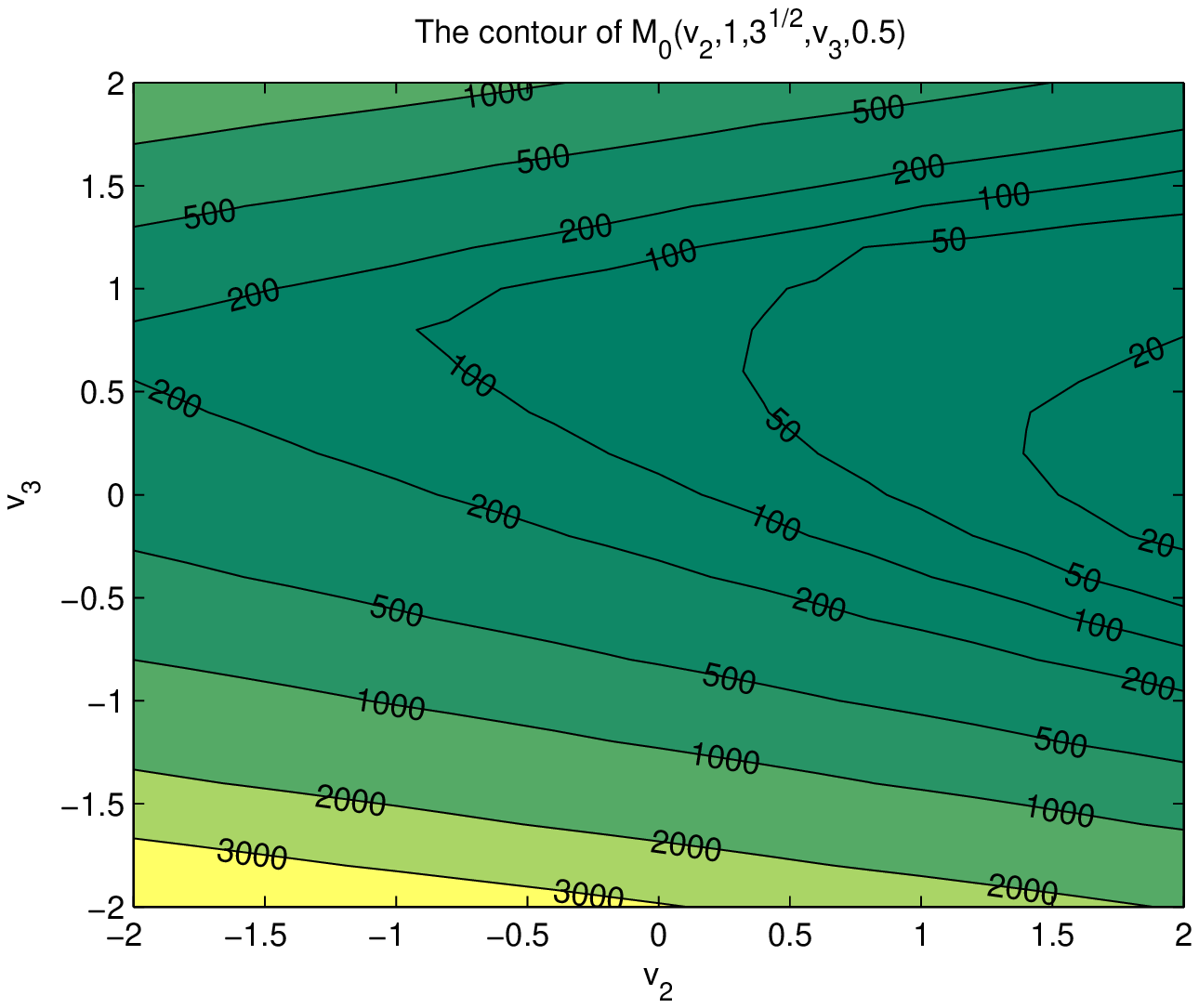} &
  \includegraphics[width=.48\textwidth]{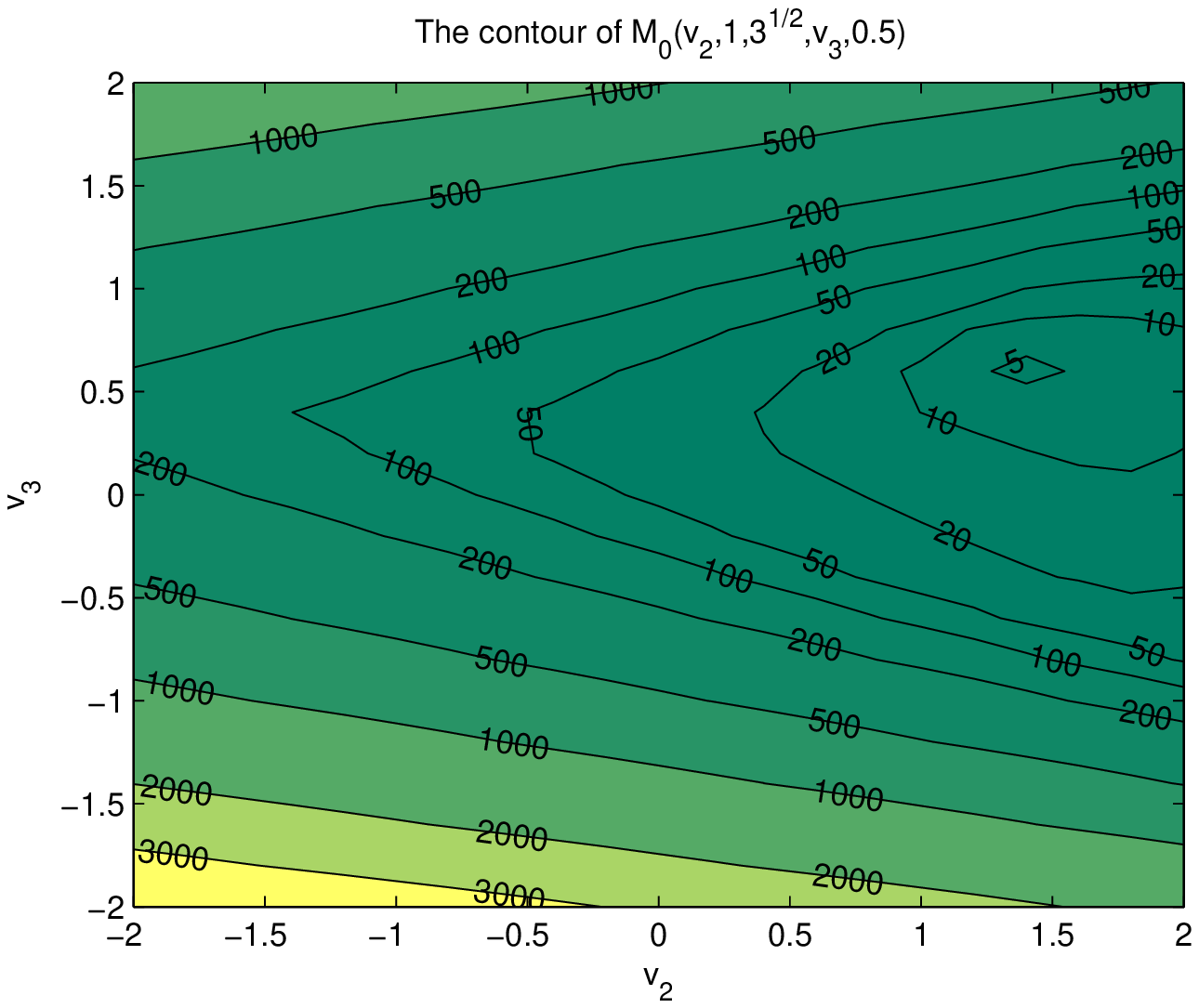} \\
  \includegraphics[width=.48\textwidth]{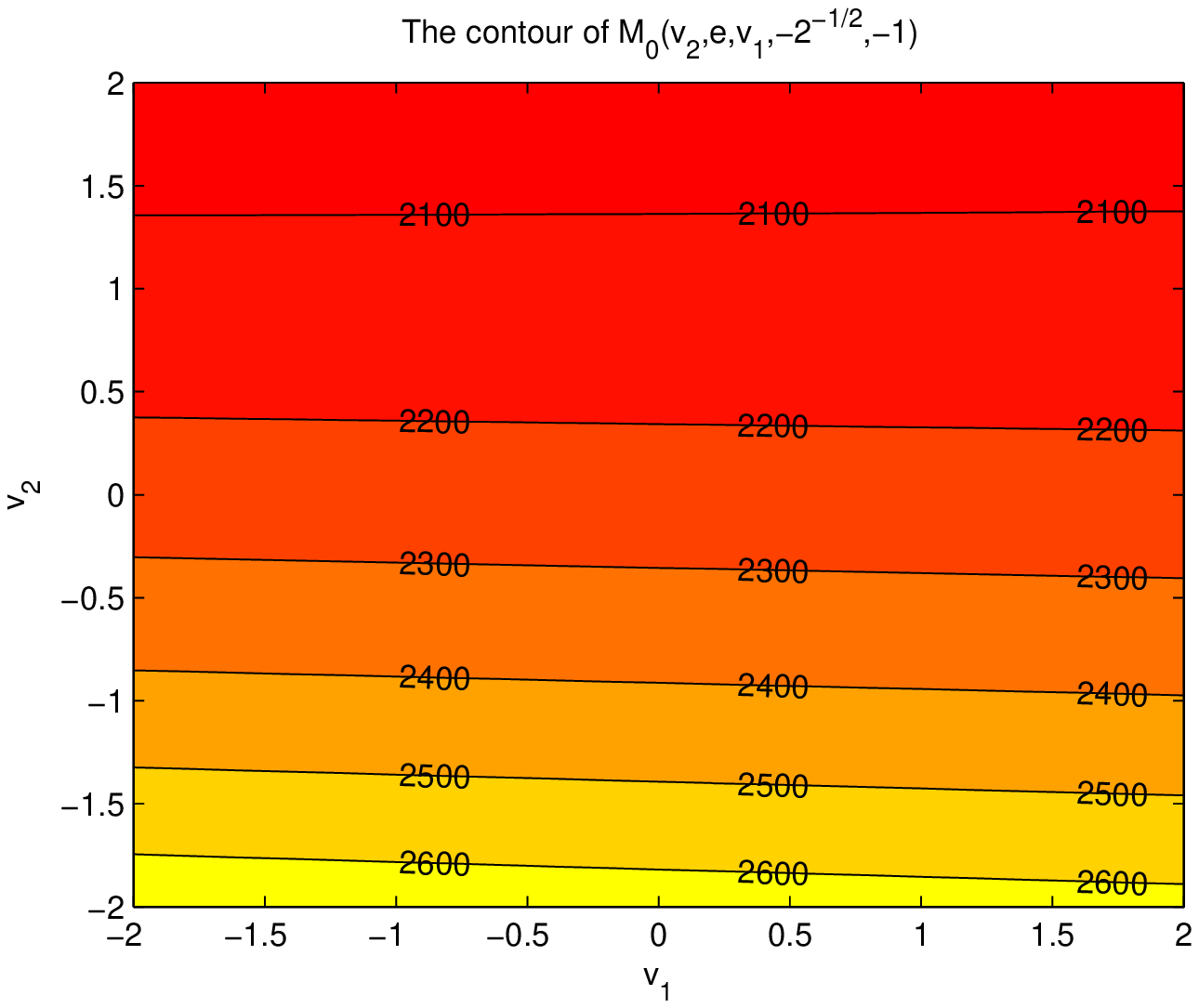} &
  \includegraphics[width=.48\textwidth]{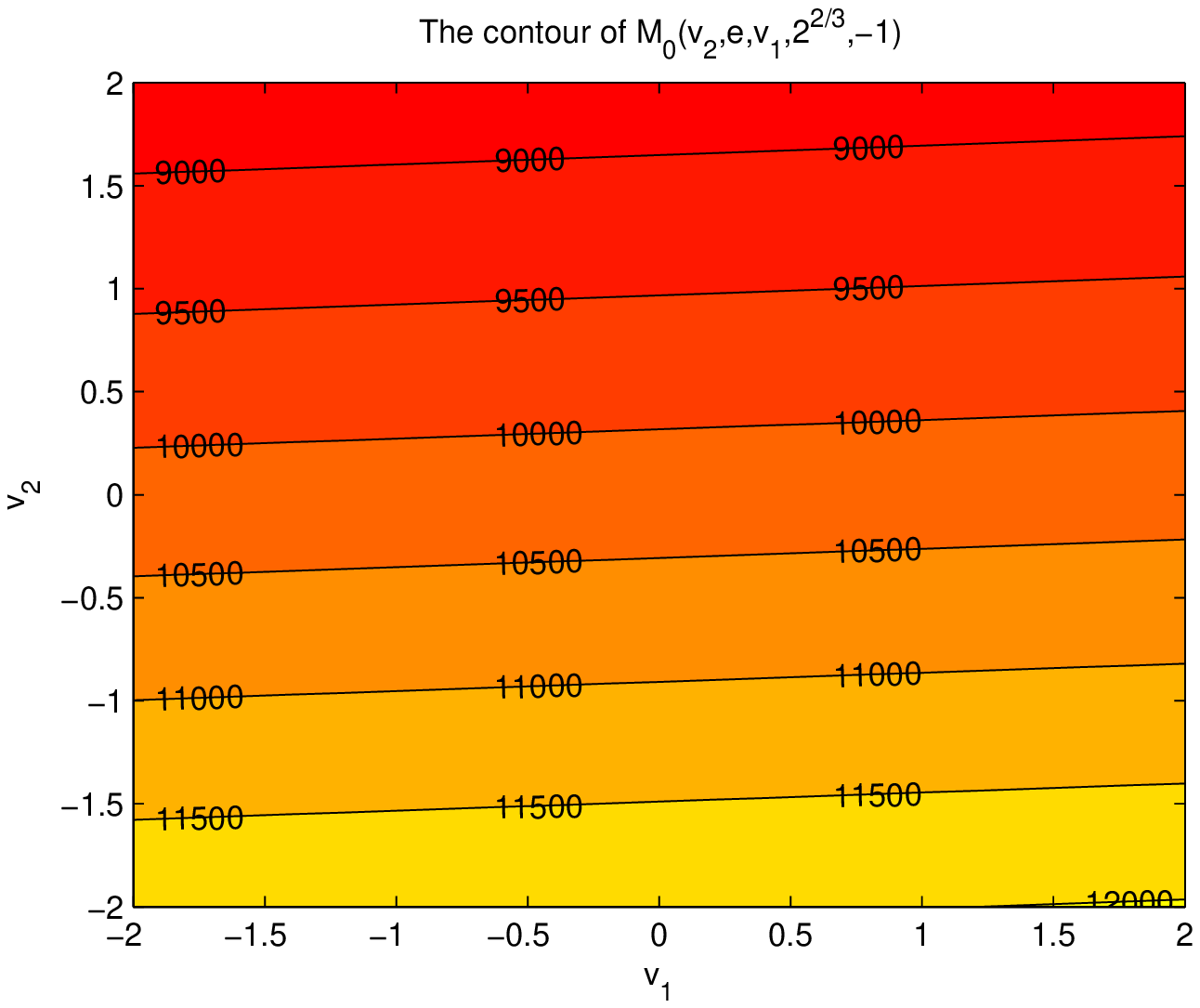}
\end{tabular}
  \caption{The contour profiles of $M_0(v_2,v_6,v_1,v_3,v_5)$ which are equivalent to $N_0(v_2,v_6,v_1,v_3,v_5)$.}\label{fig03}
\end{center}
\end{figure}

We have proved in Section 4 that some regions are PNS-free.
What about the other cases? We try to answer this problem by a numerical approach.
We use the YALMIP software with an SOS module \cite{Lo,L2} to compute $M_0(v_2,v_6,v_1,v_3,v_5)$, which is the smallest value of $v_0$ such that the fourth order four dimensional Hankel tensor $\A$ with the generating vector $(v_0,v_1,v_2,v_3,1,v_5,v_6,v_5,1,v_3,v_2,v_1,v_0)^\top$ is SOS.
Gloptipoly \cite{HLL} and SeDuMi \cite{St} are employed to compute $N_0(v_2,v_6,v_1,v_3,v_5)$, which is the smallest value of $v_0$ such that the Hankel tensor $\A$ is PSD.

\subsection{$M_0(v_2,v_6,0,0,0)$ and $N_0(v_2,v_6,0,0,0)$}

First, we focus on two elements $v_2$ and $v_6$ of generating vectors and set $v_1=v_3=v_5=0$.
By Theorem \ref{p-1}, owing to the effective domain, we have $b > -{1 \over 3}$.
We choose $v_2 = -4,-2,-1,-0.5,0,0.5,1,1.5,2,3,4$
and $v_6=-0.2,-0.1,0,0.5,1,1.5,2,4$ and compute $M_0$ and $N_0$ in these grid points respectively.
By our experiments, we found that these two functions are equivalent on all of the grid points.
Thus, no PNS tensors are detected here.
The detailed value of $M_0$ and $N_0$ are reported in Table \ref{Table-1}.


A more intuitional profile of $M_0=N_0$ is illustrated in Figure \ref{fig01}.
It is easy to see that $(v_2,v_6)=(1,1)$ is the minimizer of both $M_0$ and $N_0$
when we set $v_1=v_3=v_5=0$.

\subsection{Nonzero Odd Elements of the Generating Vectors}

We consider the case that the generating vector of a fourth order four dimensional Hankel tensor has nonzero odd elements.  According to Theorem \ref{t5}, we say that $v_5$ and $v_6$ must satisfy $\eta(v_5,v_6)<1$.
So we study them first and set $v_1=v_2=v_3=0$. We compute a plenty of grid points with different $v_5$ and $v_6$.
The function $M_0(0,v_6,0,0,v_5)$ is still equivalent to the function $N_0(0,v_6,0,0,v_5)$. That is to say, no PNS tensors are found.

The contour of $M_0(0,v_6,0,0,v_5)=M_0(0,v_6,0,0,v_5)$ is shown in Figure \ref{fig02}.
We could see that the nonlinear contour of $M_0=N_0=500$ looks like a fire balloon.

Finally, we consider all of the elements of symmetric generating vectors of fourth order four dimensional Hankel tensors. The contours of $M_0(v_2,v_6,v_1,v_3,v_5)$ and $N_0(v_2,v_6,v_1,v_3,v_5)$ for various combinations of $v_2$, $v_6$, $v_1$, $v_3$ and $v_5$ are reported in Figure \ref{fig03}.
In all of our tests, values of the function $M_0(v_2,v_6,v_1,v_3,v_5)$ in grid points are always equivalent to the corresponding values of the function $N_0(v_2,v_6,v_1,v_3,v_5)$. So,  no fourth order four dimensional PNS Hankel tensors with symmetric generating vectors are detected.

From Figures \ref{fig02} and \ref{fig03}, we could say that the second element $v_1$ of the generating vector of a Hankel tensor affect functions $M_0(v_2,v_6,v_1,v_3,v_5)$ and $N_0(v_2,v_6,v_1,v_3,v_5)$ slightly.
When we fix $v_4=1$, the middle element $v_6$ of the generating vector $\vv$ plays a more important role since it has direct impact on the effective domain.

%

\section{Final Remarks}
\hspace{4mm}

In this paper, we investigated the problem whether there exist fourth order four dimensional PNS Hankel tensors with symmetric generating vectors.   Theoretically, we proved that such PNS Hankel tensors do not exist on the segment $L = \{ (v_2, v_6, v_1, v_3, v_5)^\top = (1, 1, t, t, t)^\top: t \in [-1, 1] \}$, the cone $C = \{ (v_2, v_6, v_1, v_3, v_5)^\top =(a, b, 0, 0, 0)^\top: a \ge b \ge 1 \}$, the ray $R = \{ (v_2, v_6, v_1, v_3, v_5)^\top = (a, 0, 0, 0, 0)^\top: a \le 0 \}$ and the point  $A = (1, 0, 0, 0, 0)^\top$.  The critical value on $L$ is simply $1$.  It is interesting to note that the critical values on $C$ are a polynomial in an auxiliary parameter $\theta$ with degree four.   However, the critical values on $R$ and $A$ are irrational.   This indicate that a complete proof that fourth order four dimensional PNS Hankel tensors with symmetric generating vectors do not exist may not be easy.   However, numerical tests also indicate that such PNS Hankel tensors do not exist.   Thus, we believe that there are no fourth order four dimensional PNS Hankel tensors with symmetric generating vectors.

\vspace{3mm}
\hspace{4mm}
{\bf Acknowledgments}   We are grateful to Professor Man Kam  Kwong.  His comments helped us to improve our paper greatly.  We are also thankful to Dr. Guoyin Li for his comments.
\hspace{4mm}

%
%


\end{document}